\documentclass[12pt]{amsart}

\usepackage{latexsym, amssymb, amsmath}
\usepackage{amsmath}
 \usepackage{eucal}

    \newtheorem{rema}{Remark}[section]
    \newtheorem{propo}[rema]{Proposition}

   \newtheorem{theo}[rema]{Theorem}
 \newtheorem{conj}[rema]{Conjecture}
   
    \newtheorem{lemma}[rema]{Lemma}
    \newtheorem{corol}[rema]{Corollary}
     \newtheorem{exam}[rema]{Example}

     \newcommand{\nno}{\nonumber}
	\newcommand{\lbar}{\big\vert}
        
	\newcommand{\p}{\partial}

 \newcommand{\pf}{{\it Proof:}\hspace{2ex}}
 
 \newcommand{\epfv}{\hspace{1em}$\Box$\vspace{1em}}

\newcommand{\bC}{{\mathbb C}}

\newcommand{\bN}{{\mathbb N}}

\newcommand{\BQ}{\begin{eqnarray}}
\newcommand{\EQ}{\end{eqnarray}}
\newcommand{\BQn}{\begin{eqnarray*}}
\newcommand{\EQn}{\end{eqnarray*}}

\begin{document}
    \bibliographystyle{alpha}

\title[Reductions on Jacobian Problem in Two Variables]
{Some Reductions on Jacobian Problem in Two Variables}
    \author{Wenhua Zhao}      

\begin{abstract}
Let $f=(f_1, f_2)$ be a regular sequence of affine curves in $\bC^2$. 
Under some reduction conditions achieved by composing with some polynomial 
automorphisms of $\bC^2$,
 we show that the intersection number 
of curves $(f_i)$ in $\bC^2$ equals to the coefficient of the 
leading term $x^{n-1}$ in $g_2$, where $n=\deg f_i$ $(i=1, 2)$ and  $(g_1, g_2)$ is the
unique solution
of the equation $y{\mathcal J}(f)=g_1f_1+g_2f_2$ with $\deg g_i\leq n-1$.
So the well-known Jacobian problem is 
reduced to solving the equation above. Furthermore, by using the result above, 
we show that 
the Jacobian problem 
can also be reduced to a special family of polynomial maps. 
\end{abstract}


\maketitle


\renewcommand{\theequation}{\thesection.\arabic{equation}}
\renewcommand{\therema}{\thesection.\arabic{rema}}
\setcounter{equation}{0}
\setcounter{rema}{0}
\setcounter{section}{0}

\section{\bf Introduction}\label{S5.1}

Let $f=(f_1, f_2)$ be a pair of polynomials in two variables $(x, y)$.
Let ${\mathcal J}(f)=\frac {\p (f_1, f_2)}{\p (x, y)}$ be its Jacobian. The  
well known Jacobian Conjecture in two variables says that:  if 
${\mathcal J}(f)\equiv 1$, then the map $f: \bC^2 \to \bC^2$ is invertible and 
the inverse map $f^{-1}$ is also a polynomial map.
For the history and well known results about this conjecture, 
see \cite{BCW}  and \cite{W2}. For the two-variable case, there are 
numerous partial results. Here we just mention a few of them. 
Abhyankar \cite{A} shows that
the conjecture is equivalent to any two affine curves $(f_i)$ $(i=1, 2)$ with
the Jacobian condition having exactly one intersection point at infinity 
and gives a proof that 
any two such curves have at most two 
intersection points at infinity. He also proves the
conjecture under the condition that $k(x, y)$ is a 
Galois field extension over $k(f_1, f_2)$. Note that this is also proved by 
Markar-Limanov \cite{ML}. Nakai and Baba \cite{NB} generalize a theorem of Magnus \cite{Ma} 
and prove the conjecture
if one of $d_i=\deg f_i$ $(i=1, 2)$ is a prime, or $4$, or if  
$d_1=2p\geq d_2$ for some prime number $p$. Wright \cite{W1} proves    
that $f$ is invertible if and only if the Jacobian matrix $J(f)$ 
can be written 
as a product of elementary and diagonal matrices in $GL_2(k[X])$. 
Finally, with some help from computers, 
Moh \cite{M} shows the conjecture is true
if $d_i=\deg f_i\leq 100$ $(i=1, 2)$.

It is well known that 
the Jacobian Conjecture is equivalent to 
saying that the polynomial map $f$ is
injective when ${\mathcal J}(f)\equiv 1$ (See, for example, \cite{BCW}, \cite{R}). 
Equivalently,  if the Jacobian Conjecture is true, 
then any two affine curves $(f_i)$ $(i=1, 2)$ in $\bC^2$ 
with ${\mathcal J}(f)\equiv 1$ have (at most) one intersection 
point in $\bC^2$. 
One natural question we may ask here is whether or not
there is some  ``nice'' relationship
between 
${\mathcal J}(f)$ and the total number (counting multiplicity) of intersection 
points of the affine curves $(f_1)$ and $(f_2)$ in the affine space $\bC^2$.
Our first main result will be that, under some reduction conditions achieved by
composing with some polynomial automorphisms of $\bC^2$, 
the answer to the question above is ``Yes''.
To be more precise,
we first show that, by composing with some polynomial 
automorphisms of $\bC^2$ to 
the polynomial map $f=(f_1, f_2)$, 
we may assume that: 

$({\bf RC1})$ the leading homogeneous parts of $f_1$ and $f_2$ are both
$x^n$ for some $n\in \bN$; 

$({\bf RC2})$ all intersection points of $(f_1)$ and 
$(f_2)$ in the affine space $\bC^2$ lie
on the line $\{y=0\}$. 

From now on, we will assume $f=(f_1, f_2)$ 
satisfies the reduction conditions above 
for the rest of this section. Our first main result is 

\begin{theo}\label{MT1}
$(1)$ there is a unique polynomial solution $g=(g_1, g_2)$ for the equation
\BQ\label{E1.1}
y{\mathcal J}(f)=f_1g_1+f_2g_2
\EQ
 with $\deg g_i=n-1$ $(i=1, 2)$.

$(2)$ the total intersection number $(\mbox{counting multiplicity})$ 
of the affine curves $(f_1)$ and
$(f_2)$ in $\bC^2$ equals to the coefficient of $x^{n-1}$ of $g_2$.
\end{theo}

Since the Jacobian condition ${\mathcal J}(f)\equiv 1$  
and the total intersection number of $(f_1)$ and $(f_2)$ in $\bC^2$ can be 
preserved in our reduction procedure,   
the Jacobian conjecture in two variables is reduced to  
the problem solving the polynomial equation (\ref{E1.1}) for the polynomial
maps $f=(f_1, f_2)$ with the reduction conditions above. 
The partial solution $g_i(x, 0)$ $(i=1, 2)$ to the equation (\ref{E1.1}) 
is given in Proposition \ref{P5.3} in the case that 
all the intersection points of $(f_1)$ and $(f_2)$ in $\bC^2$ are normal
crossing. Note that the solution of the equations of the form
(\ref{E1.1}) is the so called the membership problems, which is one of the 
most important problem in computational algebra.  It has been studied  
from many different ways, see \cite{BGVY} by using Bezout Identities and residues
and  \cite{BW}, \cite{F} by using Gr\"{o}bner bases.
We hope that some results from the membership problem 
can provide some new insights to the Jacobian problem via 
Theorem \ref{MT1}.

Our second main result is 
\begin{theo}\label{MT2}
Suppose that all intersection points of $(f_1)$ and $(f_2)$ in $\bC^2$ are
normal crossing, then  the polynomial map 
$f=(f_1, f_2)$ can always be written as

\BQ\label{Form1.2}
\begin{pmatrix} f_1\\f_2
\end{pmatrix}
=
\begin{pmatrix} h_1 & -k_2\\
h_2& k_1
\end{pmatrix}
\begin{pmatrix} r(x)\\
y\lambda (x)
\end{pmatrix}
\EQ
where $\det \begin{pmatrix} h_1 & -k_2\\
h_2& k_1
\end{pmatrix}={\mathcal J}(f)$
 and $r(x)$ and $\lambda (x)$ are polynomials 
in one variable related by
\BQ
r(x)\mu(x)+r'(x)\lambda (x)= 1 \label{E1.2}\\
\deg \lambda (x) \leq \deg r(x) \label{E1.3}
\EQ
for some polynomial $\mu(x)$.
\end{theo}

From Theorem \ref{MT2}, one immediately sees that 
the Jacobian problem is reduced to 
the following 

\begin{conj}\label{Conj}
Let $f=(f_1, f_2)$ be of the form $(\ref{Form1.2})$ with the matrix 
$\begin{pmatrix} h_1 & -k_2\\
h_2& k_1
\end{pmatrix}$ being invertible. Then, for any $(r(x), \lambda (x))$
related by $(\ref{E1.2})$, $(\ref{E1.3})$ and  
 $\deg (r(x))\geq 2$, the 
Jacobian ${\mathcal J}(f)\neq c$ for any $c\in \bC^*$.
\end{conj}

The arrangement of the paper is as follows. In Section \ref{S5.2}, for the 
convenience of the readers, we fix some notation and
recall some results in the theory of residues 
and intersection numbers which will play the key roles in our later arguments.
In Section \ref{S5.3}, we first recall Noether's $AF+BG$ theorem, then
 derive some consequences which later will 
give the degree upper 
bound of the solutions of equation (\ref{E1.1}).
In Section \ref{S5.4}, we show that, by composing certain polynomial 
automorphisms of $\bC^2$, the reduction conditions $(RC1)$ and $(RC2)$ can 
be achieved. In Section \ref{S5.5}, we give the proofs for our main results 
Theorem \ref{MT1} and Theorem \ref{MT2}.

This paper is a revised version of one of the  topics in the author's Ph.D Thesis in
The University of Chicago. The author is very grateful to his adviser,
Professor Spencer Bloch for encouragement and discussions. The author
also thanks Professor Mohan Kumar, Madhav Nori and David Wright 
for personal communications on this subject.

\renewcommand{\theequation}{\thesection.\arabic{equation}}
\renewcommand{\therema}{\thesection.\arabic{rema}}
\setcounter{equation}{0}
\setcounter{rema}{0}

\section{\bf Residues and Intersection Numbers}\label{S5.2}

{\bf Notation:} 

$1)$ Let $[X_0, X_1, \cdots, X_n ]$ be the homogeneous coordinates for $\bC P^n$. 
Set $U_i=\{X_i\neq 0\}$ 
 ($i=0, 1, 2, \cdots , n$). We use $x_1, x_2, \cdots, x_n $
 to denote the Euclidean 
coordinate systems for $U_0$. We usually use small letters
 $f, g$, so on, to denote the polynomials $f(x_1, \cdots, x_n)$, 
$g(x_1, \cdots, x_n )$ in $n$ variables 
and use the corresponding capital letters to denote their homogenized polynomials 
in $X_0, X_1,\cdots, X_n$, i.e. 
$F(X_0, X_1,\cdots , X_n)=X_0^d F(\frac {X_1}{X_0}, \frac {X_2}{X_0} \cdots , 
\frac {X_n}{X_0})$,  
where $d$ is the total degree of the polynomial $f$ in $x_1, x_2, \cdots, x_n$.
 
$2)$ Let $H, F, G$ be three homogeneous polynomials in $X_i$, $(i=0, 1, 2)$. 
Suppose $F$ and $G$ intersect discretely at $p\in \bC P^2$. 
We say that {\it the restriction of $H$ at $p$
lies in the ideal generated by $F$ and $G$}, denoted by 
$H |_p \in <F, G>_p$, if the following condition is hold: 

Let $U_i$ be an affine open subset defined above for $\bC P^2$ such that 
 $p\in U_i$ for some $0\leq i\leq 2$. Here let us  assume that $i=2$.
Set $f(X_0, X_1)=F(X_0, X_1, 1)$, $g(X_0, X_1)=G(X_0, X_1, 1)$ and
$h(X_0, X_1)=H(X_0, X_1, 1)$. Then,  
as holomorphic functions near $p\in U_2$, we have $h|_p \in <f, g>_p$,
i.e. $h$ lies in the ideal generated by $f$ and $g$.

It is easy to see that the condition above does not depend on the choices of 
the affine open subset $U_i$.

A sequence $(f_1, f_2, \cdots , f_n)$, where $f_i \in {\mathcal O}_0$,
the germs of holomorphic 
functions at $0 \in \bC^n$, is said to 
be {\it regular} at $0 \in \bC^n$ if there is an open neighborhood $U$ of 
$0 \in \bC^n$ such that $0$ is the only common zeros of $f_i$ ($i=1, 2, \cdots , n$).
This is equivalent to saying that the Jacobian ${\mathcal J}(f)$ is not  identically
$0$. A sequence $(f_1, f_2, \cdots , f_n)$, where $f_i \in {\mathcal O}$
is said to be {\it regular} on $\bC^n$ if they intersect only at discrete points, or
equivalently, $f$ is regular at any point of $\bC^n$.

For a regular sequence $f=(f_1, f_2, \cdots , f_n)$ at $0\in \bC^n$ and 
a holomorphic function $h\in {\mathcal O}_0$.
Set 
\BQ
\omega (z)=\frac {h(z)}{f_1f_2\cdots f_n}dz_1\wedge dz_2\wedge \cdots \wedge dz_n
\EQ

We define
the residue of the meromorphic form at $0\in \bC^n$ to be
to be 
\BQ
{\text Res}_{\{0\}}\omega (z)=
\int_{\Gamma}\frac {h(z)}{f_1f_2\cdots f_n}dz_1\wedge dz_2\wedge \cdots \wedge dz_n
\EQ
where $\Gamma =\{z\in \bC^n : |f_1(z)|= \epsilon, |f_2(z)|= \epsilon, \cdots ,
|f_n(z)|=\epsilon \}$ for some small $\epsilon >0$.

\begin{propo}\label{P2.1}
$1)$ The residue ${\text Res}_{\{0\}}\omega $ is alternating with respect to the
permutations of $f_i$ $(i=1, 2, \cdots , n)$.

$2)$ Let $I$ be the ideal generated by $f_i$  $(i=1, 2, \cdots , n)$.
 For any $g, h \in {\mathcal O}$, set 
\BQ\label{Res2.3}
{\text Res}_{\{0\}}(g, h)=
\int_{\Gamma}\frac {g(z)h(z)}{f_1f_2\cdots f_n}dz_1\wedge dz_2\wedge \cdots \wedge dz_n
\EQ
Then, $g\in I$ if and only if  ${\text Res}_{\{0\}}(g, h)=0$ for any $h \in {\mathcal O}$. 
In other words, the bilinear pairing 
\BQ
{\text Res}_{\{0\}}: {\mathcal O}/I \times {\mathcal O}/I \rightarrow \bC
\EQ
induced by $(\ref{Res2.3})$ is  non-singular.

$3)$ {\bf (Transition Formula)} Suppose that $g=(g_1, g_2, \cdots , g_n)$ be another regular 
sequence at $0\in \bC^n$ with $g^{-1}(0)=\{0\}$. 
If $\{g_1, g_2, \cdots , g_n\}\subset I$, say 
$g_i(z)=\sum_{j=1}^n a_{i, j}(z)f_j(z)$ $(i=1, 2, \cdots , n)$, for some 
$a_{i, j}(z)\in {\mathcal O}$.
Then for any 
$h\in {\mathcal O}$, we have
\BQ
{\text Res}_{\{0\}} \left (
\frac {h(z)}{f_1f_2\cdots f_n}
dz_1\wedge dz_2\wedge \cdots \wedge dz_n\right )\nno \\  \label{TF}
  =
{\text Res}_{\{0\}}\left (
\frac {h(z) \det(A)}{g_1g_2\cdots g_n}dz_1\wedge dz_2\wedge \cdots \wedge dz_n\right )
\EQ
\end{propo}

The intersection number $D_1\cdot D_2\cdots D_n$ of the divisors $D_i=(f_i)$ 
at point $0\in \bC^n$ is
defined to be
\BQ
(D_1\cdot D_2\cdots D_n)(0) &=&
\int_{\Gamma}\frac 1{f_1f_2\cdots f_n}df_1\wedge df_2\wedge \cdots \wedge df_n \nno \\
&=& \int_{\Gamma} \frac {{\mathcal J}(f)(z)}{f_1f_2\cdots f_n}dz_1\wedge 
dz_2\wedge \cdots \wedge dz_n
\EQ

\begin{propo}\label{P2.2}
$1)$ The intersection number $(D_1\cdot D_2\cdots D_n)(0)$ 
is always a positive integer and equals to 
the degree $d$  of the holomorphic map $f$ at point $0\in \bC^n$. 
It also equals to the complex dimension of the vector space ${\mathcal O}/I$. In particular, 
for any $h\in {\mathcal O}$ with $h(0)=0$, we have $h^d \in I$.

$2)$ For any $h\in {\mathcal O}$, we have 
\BQ\label{res-deg}
{\text Res}_{\{0\}}(h, {\mathcal J}(f))={\deg}(f)h(0)
\EQ
\end{propo}

\begin{propo}\label{P2.3}
For any $h\in {\mathcal O}$ with $h(0)=0$, we have ${\mathcal J}(f)h \in I$.
\end{propo}
\pf For any $g\in {\mathcal O}$, consider
\BQ
{\text Res}_{\{0\}} (g, {\mathcal J}(f)h)
&=&{\text Res}_{\{0\}} (gh, {\mathcal J}(f)) \nno \\
&=& \deg(f)(gh)(0)\nno \\
&=& 0
\EQ
Then, by Proposition \ref{P2.1}, $2)$, we have ${\mathcal J}(f)h \in I$.
\epfv

Since the residue is defined locally, 
we can generalize it to complex manifolds. 

\begin{theo} \label{RT}{\bf (Residue Theorem)}
Let $M$ be a compact complex manifold of dimensional $n$, 
$\omega $ a meromorphic $(n, 0)$ form 
on $M$ which has only simple pole over effective divisors 
$D_i$ $(i=1, 2, \cdots , n)$. Suppose $D_i$'s intersect only at 
discrete points $v_k$ $(k=1, 2, \cdots , m)$, then
\BQ
\sum_{j=1}^m\text{Res}_{\{ v_k\}}\omega =0
\EQ
\end{theo}

\renewcommand{\theequation}{\thesection.\arabic{equation}}
\renewcommand{\therema}{\thesection.\arabic{rema}}
\setcounter{equation}{0}
\setcounter{rema}{0}

\section{\bf Noether's AF+BG Theorem and Some Consequences}\label{S5.3}

First let us  recall Noether's $AF+BG$ theorem, for the proof of this theorem, see 
\cite{GH}.
Let $H, F, G$ be three homogeneous polynomials in $X_i$, $(i=0, 1, 2)$. Let $deg(H)=m$, 
$deg(F)=k$ and $deg(G)=\ell$. 
Suppose that $a=m-k\geq 0$, $b=m-\ell \geq 0$ and $F$ and $G$ have only discrete 
intersection points. If $H |_p \in <F, G>_p$ for any 
$p\in (F)\cdot (G)$, i.e. the restriction of $H$ at any intersection point of $F$ 
and $G$ lies in the ideal generated by the restriction of $F$ and $G$. 
(See the notations fixed at the beginning of the Section \ref{S5.2}). Then there 
are homogeneous polynomials $A$ and $B$ of degree $a$ and $b$, respectively,  such 
that $H=AF+BG$. Furthermore, the pair $(A, B)$ is unique up to the following sense: 
if $(\tilde A, \tilde B)$ is another such a pair, then 
there exist a homogeneous 
polynomial $C$ such that $\tilde A=A+CG$ and $\tilde B=B-CF$. Clearly if $a<\ell$ 
or $b<k$, the pair $(A, B)$ is uniquely determined by $H, F, G$.  

Next, we derive some consequences of Noether's $AF+BG$ theorem, which will play 
crucial roles in our later argument.

\begin{propo}\label{T3.1}
Let $h, f, g$ be polynomials in $x, y$. Suppose that the affine curves 
$(f)$ and $(g)$  intersect 
only at discrete points in $\bC^2$ and $h_p\in <f, g>_p$ for any intersection point 
$p\in (f)\cdot (g)$. Then there exists a pair of polynomials 
$(a(x, y), b(x, y))$ such that 
\BQn
h(x, y)=a(x, y)f(x, y)+b(x, y)g(x, y)
\EQn
Furthermore, the pair $(a, b)$ is unique up to the similar sense as above.
\end{propo}

Note that, the proposition above is stronger than Hilbert's Nullstellensatz 
(See \cite{Ha}),
which  claims only that $h(x, y)$ is in the radical of the ideal generated 
by $f(x, y)$ and
$g(x, y)$ in general. 

\pf
We first embed the curves ${\mathcal C}_i$ for $(i=1, 2)$ into the projective 
space $\bC P^2$ by considering the homogenized polynomials $H, F, G$ of 
$h, f, g$ respectively. Observe that $F, G$ still intersect 
discretely in $\bC P^2$. 
We can choose $m\in \bN$ large enough 
such that $(X_0^m H)_p\in <F, G>_p$ for any 
$p \in (F)\cdot (G)\cdot \{X_0=0\}$ in $\bC P^2$. 
 Then apply Noether's $AF+BG$ theorem to 
the homogeneous polynomial $X_0^m H$, we have $X_0^n H=AF+BG$ for some 
homogeneous polynomials $A$ and $B$, 
then restrict to the open set $U_0\simeq \bC^2$, we get
$h=af+bg$. 
\epfv

Unlike Noether's Theorem,  Proposition \ref{T3.1}
 does not tell us much about the degrees of the polynomials  $a(x, y)$ 
and $b(x, y)$. But when $h(x, y)$ has the form 
${\mathcal J}(f)h(x, y)$, we have

\begin{propo}\label{P3.3}
Let $f_i$ $(i=1, 2)$ as above, then, for any polynomial  $h(x, y)$ which 
vanishes at all intersection points of 
$(f_1)$ and $(f_2)$ in $\bC^2$,
there exist a pair 
$(g_1, g_2)$ of polynomials  in $x, y$ such that
\BQ\label{EG3.5}
{\mathcal J}(f)h(x, y)=g_1f_1+g_2f_2
\EQ
and 
\BQ
\deg   (g_1)&\leq & \deg  (f_2)+\deg  (h)-2\\
\deg  (g_2)&\leq & \deg  (f_1)+\deg  (h)-2
\EQ
Furthermore, if $\deg ({\mathcal J}(f)h(x, y))<2\, 
\text {min}(\deg (f_1), \deg (f_2))$,
the solution $g=(g_1, g_2)$ of $(\ref{EG3.5})$ is unique.
\end{propo}

Before we give the proof of the proposition above, we need the following

\begin{lemma} \label{L3.2}
Let $f_i$ $(i=1, 2)$ as above and $F_i$  their homogenized
polynomials. Then
$\frac {\p (F_1, F_2)}{\p (X_1, X_2)} \lbar _p\in <F_1, F_2>_p$ for any
 $p\in  (F_1)\cdot (F_2)\cdot \{X_0=0\}$.
\end{lemma}
  
\pf
Let $p$ be any intersection point of $(F_1)$ and $(F_2)$ 
in $\{X_0=0\}\subset \bC P^2$. Without losing any generality, we may assume 
$p\in U_2$. Observe that for any homogeneous 
polynomial  $F(X_0, X_1, X_2)$ of degree $n$, 
we always have
\BQ
nF(X_0, X_1, X_2)&=& X_0\frac {\p}{\p X_0}F(X_0, X_1, X_2)+
X_1\frac {\p}{\p X_1}F(X_0, X_1, X_2) \nno \\
&{}& \qquad + X_2\frac {\p}{\p X_2}F(X_0, X_1, X_2)
\EQ
Therefore,

\BQ
\frac {\p}{\p X_0}F(X_0, X_1, X_2)&=&
\frac 1{X_0} \left (  nF(X_0, X_1, X_2)-
X_1\frac {\p}{\p X_1}F(X_0, X_1, X_2) \right. \nno \\
 &{}&\qquad -\left. X_2\frac {\p}{\p X_2}F(X_0, X_1, X_2) \right )
\EQ

We calculate the Jacobian $\frac {\p (F_1, F_2)}{\p (X_0, X_1)}$
as following
\BQ
&{}&\frac {\p (F_1, F_2)}{\p (X_0, X_1)}=
\begin{pmatrix}
{\frac {\p F_1}{\p X_0} } & {\frac {\p F_1}{\p X_1} }  \\
{\frac {\p F_2}{\p X_0} } & {\frac {\p F_2}{\p X_1} } \end{pmatrix} \nno \\
&=& \begin{pmatrix}
{ \frac 1{X_0} \left ( d_1F_1-
X_1\frac {\p}{\p X_1}F_1
- X_2\frac {\p}{\p X_2}F_1\right ) } & {\frac {\p F_1}{\p X_1} }  \\
{ \frac 1{X_0}\left ( d_2 F_2-
X_1\frac {\p}{\p X_1}F_2
- X_2\frac {\p}{\p X_2}F_2\right )} & {\frac {\p F_2}{\p X_1}} 
\end{pmatrix} \nno \\
&=&\frac 1{X_0}\left (d_1F_1\frac {\p F_2}{\p X_1}
-d_2F_2 \frac {\p F_1}{\p X_1}\right )
+\frac {X_2}{X_0} \frac {\p (F_1, F_2)}{\p (X_1, X_2)}
\EQ

By Proposition \ref{P2.3}, 
$X_0 \frac {\p (F_1, F_2)}{\p (X_0, X_1)}\lbar_p \in <F_1, F_2>_p$. Thus
\BQ
X_2\frac {\p (F_1, F_2)}{\p (X_1, X_2)}
&=&X_0 \frac {\p (F_1, F_2)}{\p (X_0, X_1)}-
\left (d_1F_1\frac {\p F_2}{\p X_1}
-d_2F_2 \frac {\p F_1}{\p X_1}\right ) \nno \\
&\equiv & X_0 \frac {\p (F_1, F_2)}{\p (X_0, X_1)} \quad \text{mod I}_p
\EQ
Hence $\frac {\p (F_1, F_2)}{\p (X_1, X_2)}\lbar_p \in <F_1, F_2>_p$.
\epfv

\begin{rema}\label{Rmk3.3}
If $\deg  ({\mathcal J}(f))=m$, let $J(X_0, X_1, X_2)=X_0^m{\mathcal J}(f)
(\frac {X_1}{X_0}, \frac {X_2}{X_0})$ be the homogenized polynomial of 
the Jacobian ${\mathcal J}(f)$. Then it is straightforward to check that
\BQ\label{EJj}
\frac {\p (F_1, F_2)}{\p (X_1, X_2)}
=X_0^{d_1+d_2-2-m}J(X_0, X_1, X_2)
\EQ
where $d_i=\deg  (f_i)$ for $i=1, 2$.
\end{rema}

{\it Proof of Proposition \ref{P3.3}:}  
Consider the homogeneous polynomial  
$\frac {\p (F_1, F_2)}{\p (X_1, X_2)}H$ whose restriction on $U_0$ is 
${\mathcal J}(f_1, f_2)h(x, y)$. \newline

{\it Claim:}
For any intersection point $p$  of the divisors
$(F_1)$ and $(F_2)$ in $\bC P^2$,  we have 
$\frac {\p (F_1, F_2)}{\p (X_1, X_2)}H \in <F_1, F_2>_p$. 

When $p\in U_0$, the claim follows from Proposition \ref{P2.3} 
and Remark \ref{Rmk3.3}. When 
$p \in \{ X_0=0\}$, it follows from Lemma \ref{L3.2}.

Now apply Noether's $AF+BG$ Theorem to 
$\frac {\p (F_1, F_2)}{\p (X_1, X_2)}H$, we have
\BQ
\frac {\p (F_1, F_2)}{\p (X_1, X_2)}H=
G_1F_1+G_2F_2
\EQ
for some homogeneous polynomials $G_i$ ($i=1, 2$) with 
\BQn
\deg (G_i)&=&\deg \frac {\p (F_1, F_2)}{\p (X_1, X_2)}H-\deg F_i \\
&=& \deg (F_1)+\deg (F_2)
+\deg (H)-2-\deg (F_i)
\EQn

Restrict to $U_0 \simeq \bC^2 \subset \bC P^2$, we get (\ref{EG3.5}). 

When 
$\deg ({\mathcal J}(f)h(x, y))<2\text {min}(\deg (f_1), 
\deg (f_2))$, the uniqueness of $g=(g_1, g_2)$ follows from the uniqueness 
of $G=(G_1, G_2)$ in 
Noether's $AF+BG$
theorem.
\epfv

\renewcommand{\theequation}{\thesection.\arabic{equation}}
\renewcommand{\therema}{\thesection.\arabic{rema}}
\setcounter{equation}{0}
\setcounter{rema}{0}

\section{\bf Reductions on polynomial maps}\label{S5.4}

To consider the total intersection number of a regular sequence
$f=(f_1, f_2)$  in two variables in the affine space 
$\bC^2$, we first perform the
following reductions by applying some polynomial automorphisms of $\bC^2$
with $f=(f_1, f_2)$.

Suppose the affine curves $(f_1)$ and $(f_2)$ in $\bC^2$ intersect at
discrete points $v_0, v_1, \cdots , v_N$ (without counting multiplicities).
First we can choose two 
generic lines $l_1$ and $l_2$ to form a linear basis for $\bC^2$, such that 
$v_0=l_1 \cap l_2$ and $v_i-v_j \notin l_2$ for any $i\neq j$. 
Let $(a_i, b_i)$ be the
coordinate of $v_i$ ($i=1, 2, \cdots , N$) with respect to the basis
$(l_1, l_2)$, then $v_0=(0, 0)$ and
$a_i \neq a_j$ for any $i\neq j$.

Now,  for any $m>0$, there exist a polynomial $p(x)$ such that

$1)$ $p(a_i)=b_i$ for any $i=0, 1, 2,  \cdots N$,

$2)$ $\deg  (p(x))\geq m$.

To see such a polynomial $p(x)$ always exists, we first choose and write 
$p(x)=c_{N-1} x^{m+N-1}+c_{N-2}x^{m+N-2}+\cdots +c_0 x^m$, then 
the equations $p(a_i)=b_i$ 
for $1\leq i\leq N$ give the family of the linear equations
\BQ \label{LE1}
\begin{pmatrix}
{a_1^{m+N-1}} &{a_1^{m+N-2}} &{\cdots  } &{a_1^{m}} \\
{a_2^{m+N-1}} &{a_2^{m+N-2}} &\cdots &{a_2^{m}} \\
{\vdots} &{\vdots} &{ }&{\vdots} \\
{a_N^{m+N-1}} &{a_N^{m+N-2}} &\cdots &{a_N^{m}}\end{pmatrix}
\begin{pmatrix}{c_{N-1}}\\
{c_{N-2}}\\
{\vdots }\\
{c_0}
\end{pmatrix}
=\begin{pmatrix}{b_1}\\
{b_2}\\
{\vdots }\\
{b_N}
\end{pmatrix}
\EQ

Hence (\ref{LE1}) has a unique solution. Also note that 
$p(0)=0$  since $m>0$,  therefore $p(a_i)=b_i$ for any $i=0, 1, \cdots , N$.

Observe that
 the polynomial map $u=(x, y+p(x))$ has Jacobian ${\mathcal J}(u)\equiv 1$ and is 
 invertible with the polynomial inverse $u^{-1}=(x, y-p(x))$. Therefore, the affine curves
$(f_1\circ u)$ and $(f_2\circ u)$ have same total intersection number in $\bC^2$ 
as the affine curves $(f_1)$ and $(f_2)$.
Actually, $f_1\circ u$ and $f_2\circ u$ intersect in $\bC^2$ 
only at points $(a_i, 0)$ ( $i=0, 1, 2,
\cdots , N$) which  all are  on the line $\{ y=0\}$. Another
observation is that if we choose $m$  large enough, 
the leading terms of $f_1\circ u$ and $f_2\circ u$ will depends only on $x$.
Replacing $f_1$ or $f_2$ by $f_1+f_2$ if it is necessary, we may assume that 
$f_1\circ u$ and $f_2\circ u$ have same degree. By composing certain 
linear automorphisms to $f=(f_1, f_2)$ from the left or right, which do 
not change the Jacobians, 
we may assume that the leading terms of
$f_i\circ u$ are both $x^n$.
 
From the reductions above, we see that, without losing any generality, 
for any regular sequence
$f=(f_1, f_2)$ of polynomials, composing some polynomial automorphism
to $f$ if it is necessary,
we may assume $f$ satisfies the following conditions: \newline

{\bf Reduction Conditions:}\newline

{\bf (RC1)}: $f_i=x^n+{\text {lower degree terms}}$ ( $i=1, 2$), \newline

{\bf (RC2)}: all intersection points of $(f_1)$ and 
$(f_2)$ in $\bC^2$ lie on the line $\{ y=0 \}$. \newline

\begin{rema}
From the reduction procedure above, it is easy to see that the Jacobian 
${\mathcal J}(f)$ will not be changed if we choose properly the polynomial automorphisms
composed to $f=(f_1, f_2)$. In particular, the Jacobian condition ${\mathcal J}(f)\equiv 1$
can be preserved during our reduction procedures.
\end{rema}

\renewcommand{\theequation}{\thesection.\arabic{equation}}
\renewcommand{\therema}{\thesection.\arabic{rema}}
\setcounter{equation}{0}
\setcounter{rema}{0}

\section{\bf Main Results}\label{S5.5}

From now on and for the rest of this paper, we will always 
assume $f=(f_1, f_2)$ satisfies the reduction conditions $(RC1)$ and $(RC2)$. 

Let $F_i$ $(i=1, 2)$ be the homogenized polynomial of $f_i$.
Note that, by the reduction condition $(RC1)$, the only intersection point
the curves of $(F_1)$ and $(F_2)$ in $\{ X_0=0 \}\subset \bC P^2$ is the point
$[0, 0,  1]$. We denote it by $v_\infty $.
Let $(a_i, 0)$ $(i=0, 1, 2, \cdots , N)$ be the intersection points
of $(F_1)$ and $(F_2)$ in $\bC^2\subset \bC P^2$ with $a_0=0$. 
Set $r(x)=\prod_{i=0}^{N}(x-a_i)$.   

\begin{propo}\label{P5.1}
There exist a pair of polynomials  $(k_1, k_2)$ and a unique pair of polynomials
$(g_1, g_2)$ such that
\BQ
{\mathcal J}(f) y   &=& g_1f_1+g_2f_2 \label{GE5.1}\\
{\mathcal J}(f) r(x)&=& k_1f_1+k_2f_2 \label{KE5.1}
\EQ
with $deg (g_i)\leq n-1$ and $deg (k_i)\leq n+N-1$ $(i=1, 2)$.
\end{propo}

\pf This is a direct consequence of Proposition \ref{P3.3} for the 
polynomials $h(x, y)=y$ and $h(x, y)=r(x)$, respectively. 
\epfv

\begin{theo} \label{T5.2}
With the same notation as above, 
the coefficient of the term $x^{n-1}$ of $g_2$ equals to
the intersection number of $f_1$ and $f_2$ in $\bC^2$.
\end{theo}

Note that Proposition \ref{P5.1} and Theorem \ref{T5.2} imply Theorem \ref{MT1}, 
the first main result stated in Section \ref{S5.1}.

\pf Let $m=\deg  ({\mathcal J}(f))$ and 
$J=X_0^m{\mathcal J} (f)(\frac {X_1}{X_0}, \frac {X_2}{X_0})$. 
Note that under the reduction conditions 
$(RC1)$ and $(RC2)$, we have $m<2n-2$.

Consider the meromorphic $(2, 0)$ form 
$\omega$ on $\bC P^2$ which is defined as following:

\BQ\label{omega} 
\omega (x, y)=\frac {{\mathcal J}(f)}{f_1f_2}  dx \wedge dy 
\quad \text {on} \quad U_0 
\EQ

\BQn
\omega (x_0, x)=\frac {x_0^{2n-3-m}J(x_0, x)}
{F_1(x_0, x)F_2(x_0, x)}  dx_0 \wedge dx \quad \text {on} \quad U_2
\EQn

\BQn
 \omega (x_0, y)=-\frac {x_0^{2n-3-m}J(x_0, y)}
{F_1(x_0, y)F_2(x_0, y)} dx_0 \wedge dy \quad \text {on} \quad U_1
\EQn
where $(x, y)$, $(x_0, x)$ and $(x_0, y)$ are the Euclidean coordinates
for $U_0$, $U_2$ and $U_1$, respectively.

It is easy to check that $\omega$ is well defined $(2, 0)$ form 
on $\bC P^2$ and has pole only at the effective divisors $(F_1)$ and $(F_2)$. 
Then by the Residue Theorem \ref{RT}, we have
\BQ \label{PP5.4}
{Res}_{v_\infty} \omega =-
\sum_{i=0}^{i=N}\text {Res}_{v_i} \omega 
\EQ
which, by  (\ref{omega}), 
is the negative of the intersection number
 of $(f_1)$ and $(f_2)$
in $\bC^2$.

We can calculate the residue of $\omega$ at $v_\infty$ as following:
Note that from Remark \ref{Rmk3.3} and (\ref{GE5.1}), we have
\BQn
X_0^{2n-2-m}J Y=G_1F_1+G_2F_2=
(G_1+G_2)F_1+G_2(F_2-F_1)
\EQn
Note that
\BQ
F_1&=&X_1^n+X_0B_1 \\
F_2&=&X_1^n+X_0B_2
\EQ
for some homogeneous polynomials $B_i$. So
$F_2-F_1=X_0(B_2-B_1)$ is divisible by $X_0$, hence so is 
$G_1+G_2$. Let $B=X_0^{-1}(F_2-F_1)$, then
\BQ
X_0^{2n-3-m} JY=X_0^{-1}(G_1+G_2)F_1+G_2B
\EQ
Restrict  to $U_2$, we have
\BQn
x_0^{2n-3-m} J(x_0, x)=x_0^{-1}(G_1+G_2)(x_0, x)F_1(x_0, x)+
G_2(x_0, x)B(x_0, x)
\EQn
By Proposition \ref{P2.1}, the Transition Formula (\ref{TF}), we have
\BQn
\text {Res}_{v_\infty}\omega &=&\int_{v_\infty}
 \frac {x_0^{2n-3-m} J(x_0, x)}
{F_1(x_0, x)F_2(x_0, x)}dx_0\wedge dx\\
&=&
\int_{v_\infty}\frac {G_2 B}{F_1(x_0, x)F_2(x_0, x)}dx_0\wedge dx\\
&=&
-\int_{v_\infty}\frac {G_2}{x^nx_0}dx_0\wedge dx\\
&=& -\text{the coefficient of $x^{n-1}$ of $G_2$}.
\EQn
Then by (\ref{PP5.4}), we are done.
\epfv

It is interesting and probably a little surprising to see that
the total intersection number of $(f_1)$ and $(f_2)$ in $\bC^2$ and
the Jacobian ${\mathcal J}(f)$ are related in the algebraic way provided by
Proposition \ref{P5.1} and Theorem \ref{T5.2}. Considering the Jacobian problem, 
we immediately have

\begin{corol}\label{CC1}
The Jacobian Conjecture for two variables is equivalent to 
the following statement:
Suppose that $f=(f_1, f_2)$ satisfies the Reduction Conditions 
$(RC1)$, $(RC2)$ and ${\mathcal J}(f)\equiv 1$. Let $g=(g_1, g_2)$ be 
the unique solution of the equation 
\BQ\label{E5.8}
y=g_1f_1+g_2f_2
\EQ
with $\deg (g_i)\leq n-1$. Then the coefficient of $x^{n-1}$
of $g_2$ equals $1$.
\end{corol}

Unfortunately, the equation (\ref{GE5.1}) is not 
quite easy to solve in general,
even though it is a 
 linear equation 
and has a unique solution $g=(g_1, g_2)$ with the degree condition 
$\deg g_i\leq n-1$ $(i=1, 2)$. One question, which we think, might be interesting is
to look more closely at the algorithm using Gr\"{o}bner bases and to see if we can get 
more insight 
to solution of the equation (\ref{GE5.1}) or (\ref{E5.8}).
 
In the next proposition,  we give the partial solution 
 $g_i(x, 0)$ 
of $g_i$ in terms of $f_i(x, 0)$ ($i=1, 2$) under 
 the condition 
that $(f_1)$ and $(f_2)$  have only transversal
intersection points in $\bC^2$.

Write $f_1=\sum_{i=0}^{n-1} a_i(x)y^i$,
$f_2=\sum_{i=0}^{n-1} b_i(x)y^i$,
$g_1=\sum_{i=0}^{n-1} c_i(x)y^i$ and 
$g_2=\sum_{i=0}^{n-1} d_i(x)y^i$, then we have

\begin{propo}\label{P5.3}
Suppose  
${\mathcal J}(f)(v_i) \neq 0$ for any $i=0, 1,2, \cdots , N$.
Then,
\BQ
 c_0(x)&=&-r'(x)\frac {b_0(x)}{r(x)} \label{E5.9}\\
 d_0(x)&=& r'(x)\frac {a_0(x)}{r(x)}\label{E5.10}
\EQ
where $r(x)=\Pi_{i=0}^{N}(x-a_i)$ as before and $r'(x)=\frac {dr}{dx}(x)$.
\end{propo}
\pf  
From the equations (\ref{GE5.1}) and (\ref{KE5.1}), let $y=0$, we get

\BQ
0 &=&c_0(x)a_0(x)+d_0(x)b_0(x) \label{E5.12}\\ \label{E5.13}
{\mathcal J}(f)(x, 0) r(x)&=&k_1(x, 0)a_0(x)+k_2(x, 0)b_0(x)
\EQ

Note that $r(x)$ is the greatest common divisor of $a_0(x)$ and $b_0(x)$,
 so $a_0(x)/r(x)$ and $b_0(x)/r(x)$ 
are coprime to each
other. Dividing $r(x)$ from the both sides of (\ref{E5.12}), we get
$ 0 =c_0(x)(a_0(x)/r(x))+d_0(x)(b_0(x)/r(x))$. Therefore there exists 
a polynomial $\eta (x)$ such that
$c_0(x)=-\eta (x)(b_0(x)/r(x))$, 
$d_0(x)=\eta (x)(a_0(x)/r(x))$  and $\deg (\eta (x))\leq N$.
It is easy to check that 
\BQ\label{E5.14}
k_1(x, 0)d_0(x)-k_2(x, 0)a_0(x)=
{\mathcal J}(f)(x, 0)\eta (x)
\EQ

Now apply the Transition Formula (\ref{TF}) to 
the equations (\ref{GE5.1}) and (\ref{KE5.1}), we get
\BQn \label{PE15}
1&=&\text{Res}_{\{ v_i\} } (\frac {{\mathcal J}(f)}{f_1f_2} dx\wedge dy) \\
&=&\text{Res}_{\{ v_i\}}   (\frac {{\mathcal J}(f)(g_1k_2-g_2k_1)}
{({\mathcal J}(f)y)({\mathcal J}(f)r(x))} dx\wedge dy) \nno \\
&=& \text{Res}_{\{ v_i\}}  (\frac { (g_1k_2-g_2k_1){\mathcal J}(f)^{-1} }
{yr(x)} dx\wedge dy) \nno \\
&=& \text{Res}_{\{ v_i\}} (\frac { (g_1k_2-g_2k_1)(x, 0){\mathcal J}(f)^{-1}(x, 0) }
{yr(x)} dx\wedge dy) \nno \\
&=& -\text{Res}_{\{ v_i\}} (\frac {(k_1(x, 0)d_0(x)-k_2(x, 0)a_0(x))
{\mathcal J}(f)^{-1}(x, 0)}
{yr(x)} dx\wedge dy) \nno \\
&=& -\text{Res}_{\{ v_i\}} (\frac {\eta (x)}
{yr(x)} dx\wedge dy) \nno \\ 
&=& \frac {\eta (a_i)}{r'(a_i)} 
\EQn

Hence,  we have 
$\eta (a_i)=r'(a_i)$ for any $i=0, 1, 2, \cdots , N$.
Since $\deg ( \eta (x) )\leq \deg (r'(x))$. 
we have $\eta (x)=r'(x)$.  
\epfv

Under the conditions above, we can choose the following special solution 
$k=(k_1, k_2)$ for 
the equation (\ref{GE5.1}) as follows: 

Set
\BQ
k_1(x, y)&=&\frac 1{y}(r(x)g_1(x, y) + r'(x)f_2(x, y)) \label{5E1}\\
k_2(x, y)&=&\frac 1{y}(r(x)g_2(x, y) - r'(x)f_1(x, y))\label{5E2}
\EQ 
Note that from equations (\ref{E5.9}) and (\ref{E5.10}), it is easy to check that
$k_i$ ($i=1, 2$) defined above are polynomials.

From (\ref{5E1}) $\times f_1 +$ (\ref{5E2})$ \times f_2$, we get
\BQ\label{5E3}
{\mathcal J}(f)r(x)=k_1(x, y)f_1(x, y) +k_2(x, y)f_2(x, y)
\EQ

From (\ref{5E1}) $\times g_2 - $ (\ref{5E2}) 
$ \times g_1$, we get 

\BQ\label{5E4}
{\mathcal J}(f)r'(x)=k_1(x, y)g_2(x, y)-k_2(x, y) g_1(x, y)
\EQ

Since that $r(x)$ and $r'(x)$ are coprime to each other, there exist 
polynomials $\lambda (x)$ and $\mu (x)$ such that

$1)$ 
\BQ
r(x)\mu (x)+r'(x)\lambda (x)&=& 1
\EQ

$2)$ $\deg  \lambda (x)\leq N$ and $\deg  \mu (x)\leq N-1$.

Set 
\BQ
h_1(x, y)&=&\mu (x)f_1(x, y)+\lambda (x)g_2(x, y) \\ \label{PEH}
h_2(x, y)&=&\mu(x)f_2(x, y) -\lambda (x)g_1(x, y)
\EQ
Then, from the equations (\ref{5E1})-(\ref{PEH}), we have
\BQ
&{}& \det \begin{pmatrix} h_1 & -k_2\\ h_2& k_1 \end{pmatrix} 
={\mathcal J}(f)
\label{PE5} \\
&{}& \begin{pmatrix} f_1\\f_2 \end{pmatrix}
=
\begin{pmatrix} h_1 & -k_2\\ h_2& k_1 \end{pmatrix} 
\begin{pmatrix} r(x)\\ y\lambda (x)\end{pmatrix} \label{PE10}\\
&{}& \begin{pmatrix} g_1\\g_2 \end{pmatrix}
=
\begin{pmatrix} -h_2 & k_1\\ h_1& k_2 \end{pmatrix}
\begin{pmatrix} r'(x)\\ y\mu (x) \end{pmatrix}\label{PE12}
\EQ

In particular,  we have proved Theorem \ref{MT2}, the second main result stated 
in Section \ref{S5.1}.

\begin{exam}
Consider $f=(f_1, f_2)$, where
\BQn
f_1(x, y)&=&x+y+x^n\\
f_2(x, y)&=&y+x^n
\EQn
Note that ${\mathcal J}(f)\equiv 1$ and
$f$ is a polynomial automorphism of $\bC^2$  with the inverse 
$f^{-1}=(f_1^{-1}, f_2^{-1})$, where
\BQn
f_1^{-1}(x, y)&=&x-y\\
f_2^{-1}(x, y)&=&y-(x-y)^n
\EQn
Hence,  $(f_1)$ and $(f_2)$ intersect only at  $0\in \bC^2$ 
with multiplicity $1$ and in this case, we have
$r(x)=x$, $\mu (x)=0$ and $\lambda (x)=1$. 

It is easy to check that the unique solution $g=(g_1, g_2)$ of equation $(\ref{E5.8})$
is given by
\BQn
g_1(x, y)&=&-x^{n-1}\\
g_2(x, y)&=&x^{n-1}+1
\EQn
Hence the 
coefficient of $x^{n-1}$ of $g_2$ is same as  the total intersection number 
of $(f_1)$ and $(f_2)$ in $\bC^2$ which is $1$.  

We choose the matrix
\BQn
\begin{pmatrix} h_1 & -k_2\\
h_2 & k_1
\end{pmatrix}
=\begin{pmatrix} 1+x^{n-1} & 1\\
x^{n-1} & 1
\end{pmatrix}
\EQn
Then we have $\det \begin{pmatrix} h_1 & -k_2\\
h_2 & k_1
\end{pmatrix}={\mathcal J}(f)\equiv 1$ and  

\BQn
\begin{pmatrix} f_1\\f_2
\end{pmatrix}
&=&
\begin{pmatrix} 1+x^{n-1} & 1\\
x^{n-1} & 1
\end{pmatrix}
\begin{pmatrix} x \\ y \end{pmatrix}\\
\begin{pmatrix} g_1\\g_2
\end{pmatrix}
&=&
\begin{pmatrix} -x^{n-1} & 1\\
1+x^{n-1} & -1
\end{pmatrix}
\begin{pmatrix} 1\\ 0 \end{pmatrix}
\EQn
which are the equations $(\ref{PE10})$ and $(\ref{PE12})$, respectively, in this case.

\end{exam}

{\small \sc Department of Mathematics, Washington University in St. Louis,
St. Louis, MO 63130 }

{\em E-mail}: zhao@math.wustl.edu

\end{document}